\documentclass[11pt,fleqn]{article}
\usepackage{amsmath, amssymb, graphicx}
\usepackage{latexsym}
\usepackage{xcolor}

\begin{document}
\newcommand{\Bo}{\rule{2mm}{3mm}}
\newcommand{\be}{\begin{equation}}
\newcommand{\ee}{\end{equation}}
\newcommand{\mathsym}[1]{{}}
\newcommand{\unicode}{{}}
\renewcommand{\theequation}{\thesection.\arabic{equation}}


\title{\Large{\bf Order of Convexity of Integral Transforms and Duality}}
{\author{Sarika Verma $^\dagger$,  Sushma Gupta and Sukhjit Singh \\
 Sant Longowal Institute of Engineering and Technology,\\
 Longowal-148106 (Punjab)-INDIA.\\
\rm{ e-mail: $^\dagger$sarika.16984@gmail.com }}}

\date{}
\maketitle

\begin{abstract}
\noindent Recently, Ali et al \cite{rmali1} defined the class $\mathcal{W}_{\beta}(\alpha, \gamma)$ consisting of functions $f$ which satisfy
$$\Re e^{i\phi}\left((1-\alpha+2\gamma)\frac{f(z)}{z}+(\alpha-2\gamma)f'(z)+\gamma zf''(z)-\beta\right)>0,$$
for all $z\in E=\left\{z : |z|<1\right\}$ and for $\alpha, \gamma\geq0$ and $\beta<1$, $\phi\in \mathbb{R}$ (the set of reals). For $f\in{\mathcal{W}_{\beta}(\alpha, \gamma)}$, they discussed the convexity of the integral transform $$V_{\lambda}(f)(z):=\int_{0}^{1}\lambda(t)\frac{f(tz)}{t}dt,$$
where $\lambda$ is a non-negative real-valued integrable function satisfying the condition $\displaystyle\int_{0}^{1}\lambda(t)dt=1$. The aim of present paper is to find conditions on $\lambda(t)$ such that $V_{\lambda}(f)$ is convex of order $\delta$ ($0\leq\delta\leq1/2$) whenever $f\in{\mathcal{W}}_{\beta}(\alpha, \gamma)$. As applications, we study various choices of $\lambda(t)$, related to classical integral transforms.
\end{abstract}




\section{Introduction}
\setcounter{equation}{0}

Let ${\mathcal{A}}$ denote the class of analytic functions $f$ defined in the open unit disc $E=\{z : |z|<1 \}$ with the normalization $f(0)=f'(0)-1=0$. Let ${\mathcal{A}}_{0}=\left\{g: \, g(z)=f(z)/z, \, f\in {\mathcal{A}}\right\}$. Let $S$ be the subclass of $\mathcal{A}$ consisting of univalent functions in $E$. A function $f\in S$ is said to be starlike or convex, if f maps $E$ conformally onto the domains, respectively, starlike with respect to the origin and convex. The generalization of these two classes are given by the following analytic characterizations :
$$ S^{\ast}(\beta)=\left\{f\in{\mathcal{A}} : \Re\left(\frac{zf'(z)}{f(z)}\right)>\beta, \hspace{0.5cm}  0\leq\beta<1\right\}$$
$$ K(\beta)=\left\{f\in{\mathcal{A}} : \Re\left(1+\frac{zf''(z)}{f'(z)}\right)>\beta, \hspace{0.5cm}  0\leq\beta<1\right\}.$$
For $\beta=0$, we usually set $S^{\ast}(0)=S^{\ast}$ and $K(0)=K$.\\

For two functions $f(z)=z+a_{2}z^{2}+a_{3}z^{3}+\cdots$ and $g(z)=z+b_{2}z^{2}+b_{3}z^{3}+\cdots$ in $\mathcal{A}$, their Hadamard product (or convolution) is the function $f \ast g$ defined by
$$(f \ast g)(z)=z+ \sum_{n=2}^\infty a_nb_nz^n.$$

For $f\in {\mathcal{A}}$, Fournier and Ruscheweyh \cite{fournier} introduced the operator
\begin{equation}
F(z)=V_{\lambda}(f)(z):=\int_{0}^{1}\lambda(t)\frac{f(tz)}{t}dt,
\end{equation}
where $\lambda$ is a non-negative real-valued integrable function satisfying the condition $\displaystyle\int_{0}^{1}\lambda(t)dt=1$. This operator
contains some of the well-known operators such as Libera, Bernardi and Komatu as its special cases. This operator has been studied by a number of authors for various choices of $\lambda(t)$ (for example see \cite{rmali}, \cite{bala}, \cite{bala-dj}, \cite{fournier}). Fournier and Ruscheweyh \cite{fournier} applied the duality theory (\cite{ruscheweyh,ruscheweyh1}) to prove the starlikeness of the linear integral transform $V_{\lambda}(f)$ when $f$ varies in the class
$${\mathcal{P}}(\beta):=\left\{f\in {\mathcal{A}}:\exists\phi\in{\mathbb{R}} | \Re e^{i\phi}\left(f'(z)-\beta\right)>0, \, \, z\in E\right\}.$$

In 1995, Ali and Singh \cite{alisingh} discussed the convexity properties of the integral transform (1.1) for functions $f$ in the class ${\mathcal{P}}(\beta)$. In 2002, Choi et al. \cite{choikim} investigated convexity properties of the integral transform (1.1) for functions $f$ in the class $${\mathcal{P}}_{\gamma}(\beta):=\left\{f\in {\mathcal{A}}:\exists\phi\in{\mathbb{R}} | \Re e^{i\phi}\left((1-\gamma)\frac{f(z)}{z}+\gamma f'(z)-\beta\right)>0, \, \, z\in E\right\}.$$
It is evident that the class ${\mathcal{P}}_{\gamma}(\beta)$ is closely related to the class ${\mathcal{R}}_{\gamma}(\beta)$ defined by $${\mathcal{R}}_{\gamma}(\beta):=\left\{f\in {\mathcal{A}}:\exists\phi\in{\mathbb{R}} | \Re e^{i\phi}\left(f'(z)+\gamma zf''(z)-\beta\right)>0, \, z\in E\right\}.$$ Clearly, $f\in{\mathcal{R}}_{\gamma}(\beta)$ if and only if $zf'$ belongs to ${\mathcal{P}}_{\gamma}(\beta)$.\\

\indent In a very recent paper, R.M.ali et al \cite{rmali1} discussed the convexity of the integral transform (1.1) for the functions $f$ in a more general class $\mathcal{W}_{\beta}(\alpha, \gamma)$
\begin{equation}
\left\{f\in {\mathcal{A}}:\exists\phi\in{\mathbb{R}} | \Re e^{i\phi}\left((1-\alpha+2\gamma)\frac{f(z)}{z}+(\alpha-2\gamma)f'(z)+\gamma zf''(z)-\beta\right)>0, \, z\in E\right\}.
\end{equation}
\noindent Note that ${\mathcal{W}_{\beta}(1, 0)}\equiv{\mathcal{P}}(\beta)$, ${\mathcal{W}_{\beta}(\alpha, 0)}\equiv{\mathcal{P}}_{\alpha}(\beta)$
and ${\mathcal{W}_{\beta}(1+2\gamma, \gamma)}\equiv{\mathcal{R}}_{\gamma}(\beta)$.\\

\indent In the present paper, we shall mainly tackle the problem of finding a sharp estimate of the parameter $\beta$ that ensures $V_{\lambda}(f)$ to be convex of order $\delta$ for $f\in{\mathcal{W}_{\beta}(\alpha, \gamma)}$. To prove our result, we shall need the duality theory for convolutions, so we include here some basic concepts and results from this theory. For a subset $\mathcal{B}\subset{\mathcal{A}}_{0}$, we define
$${\mathcal{B}}^{*}=\left\{g\in{\mathcal{A}}_{0}: \, (f\ast g)(z)\neq0, \, z\in E, \, \text{for all} \, f\in{\mathcal{B}}.\right\}$$
The set ${\mathcal{B}}^{*}$ is called the dual of ${\mathcal{B}}$. Further, the second dual of ${\mathcal{B}}$ is defined as ${\mathcal{B}}^{**}=({\mathcal{B}}^{*})^{*}$. We state below a fundamental result.\\

\textbf{Theorem 1.1.} Let
$${\mathcal{B}}=\left\{\beta+(1-\beta)\left(\frac{1+xz}{1+yz}\right): \, |x|=|y|=1\right\}, \, \, \beta\in{\mathbb{R}}, \, \beta\neq1.$$
Then, we have\\

\indent(1) ${\mathcal{B}}^{**}=\left\{g\in{\mathcal{A}}_{0}: \, \exists \phi\in{\mathbb{R}} \, \text{such that} \, \Re\{e^{i\phi}(g(z)-\beta)\}>0, \, \, z\in E\right\}$.\\

\indent(2) If $\Gamma_{1}$ and $\Gamma_{2}$ are two continuous linear functionals on $\mathcal{B}$ with $0{\not\in}\Gamma_{2}$, then for every $g\in{\mathcal{B}}^{**}$ we can find ${v}\in{\mathcal{B}}$ such that
$$\frac{\Gamma_{1}(g)}{\Gamma_{2}(g)}=\frac{\Gamma_{1}(v)}{\Gamma_{2}(v)}.$$

 The basic reference to this theory is the book by Ruscheweyh \cite{ruscheweyh} (see also \cite{ruscheweyh1}).


\section{Preliminaries}
\setcounter{equation}{0}

We follow the notations used in \cite{rmali}. Let $\mu\geq0$ and $\nu\geq0$ satisfy
\begin{equation}
\mu+\nu=\alpha-\gamma \, \, \, \, \text{and} \, \, \, \, \mu\nu=\gamma.
\end{equation}
When $\gamma=0$, then $\mu$ is chosen to be 0, in which case, $\nu=\alpha\geq0$. When $\alpha=1+2\gamma$, (2.1) yields $\mu+\nu=1+\gamma=1+\mu\nu$, or $(\mu-1)(1-\nu)=0$.\\

\noindent (i)  For $\gamma>0$, then choosing $\mu=1$ gives $\nu=\gamma$.\\
(ii) For $\gamma=0$, then $\mu=0$ and $\nu=\alpha=1$.\\

Whenever the particular case $\alpha=1+2\gamma$ will be considered, the values of $\mu$ and $\nu$ for $\gamma>0$ will be taken as $\mu=1$ and $\nu=\gamma$ respectively, while $\mu=0$ and $\nu=1=\alpha$ in the case when $\gamma=0$.\\

\indent Next we introduce two auxiliary functions. Let
\begin{equation}
{\phi}_{\mu,\nu}(z)=1+\sum_{n=1}^{\infty}\frac{(n\nu+1)(n\mu+1)}{n+1}z^{n},
\end{equation}
and
\begin{eqnarray} {\psi}_{\mu,\nu}(z) & = & {\phi}_{\mu,\nu}^{-1}(z)=1+\sum_{n=1}^{\infty}\frac{n+1}{(n\nu+1)(n\mu+1)}z^{n} \nonumber\\
& = & \int_{0}^{1}\int_{0}^{1}\frac{dsdt}{(1-t^{\nu}s^{\mu}z)^2}. \end{eqnarray}

Here ${\phi}_{\mu,\nu}^{-1}$ denotes the convolution inverse of ${\phi}_{\mu,\nu}$ such that ${\phi}_{\mu,\nu}\ast{\phi}_{\mu,\nu}^{-1}=z/(1-z)$. If $\gamma=0$, then $\mu=0$, $\nu=\alpha$, and it is clear that
$${\psi}_{0,\alpha}(z)=1+\sum_{n=1}^{\infty}\frac{n+1}{n\alpha+1}z^{n}=\int_{0}^{1}\frac{dt}{(1-t^{\alpha}z)^2}.$$
If $\gamma>0$, then $\nu>0$, $\mu>0$, and making the change of variables $u=t^{\nu}$, $v=s^{\mu}$ results in
$${\psi}_{\mu,\nu}(z)=\frac{1}{\mu\nu}\int_{0}^{1}\int_{0}^{1}\frac{u^{1/{\nu}-1}v^{1/{\mu}-1}}{(1-uvz)^2}dudv.$$
Thus the function ${\psi}_{\mu,\nu}$ can be written as\\
\begin{equation}{\psi}_{\mu,\nu}(z)=\left\{
                                      \begin{array}{ll}
                                        \frac{1}{\mu\nu}\int_{0}^{1}\int_{0}^{1}\frac{u^{1/{\nu}-1}v^{1/{\mu}-1}}{(1-uvz)^2}dudv, & \hbox{$\gamma>0$;} \\
                                        \int_{0}^{1}\frac{dt}{(1-t^{\alpha}z)^2}, & \hbox{$\gamma=0$, $\alpha>0$.}
                                      \end{array}
                                    \right.
\end{equation}

Let $q$ be the solution of the initial value problem
\begin{equation}\frac{d}{dt}\left(t^{1/{\nu}}q(t)\right)=\left\{
                                      \begin{array}{ll}
                                        \frac{1}{\mu\nu}t^{1/{\nu}-1}\int_{0}^{1}\frac{(1-\delta)-(1+\delta)st}{(1-\delta)(1+st)^{3}}s^{1/{\mu}-1}ds, & \hbox{$\gamma>0$,} \\
                                        \frac{1}{\alpha}\frac{(1-\delta)-(1+\delta)t}{(1-\delta)(1+t)^{3}}t^{1/{\alpha}-1}, & \hbox{$\gamma=0$, $\alpha>0$,}
                                      \end{array}
                                    \right.
\end{equation}
satisfying $q(0)=1$.\\

\noindent Solving the differential equation (2.5), we have
\begin{equation}
q(t)=\frac{1}{\mu\nu}\int_{0}^{1}\int_{0}^{1}\frac{(1-\delta)-(1+\delta)swt}{(1-\delta)(1+swt)^{3}}s^{1/{\mu}-1}w^{1/{\nu}-1}dsdw.
\end{equation}
In particular,
\begin{equation}
q_{\alpha}(t)=\frac{1}{\alpha}\int_{0}^{1}\frac{(1-\delta)-(1+\delta)st}{(1-\delta)(1+st)^{3}}s^{1/{\alpha}-1}ds, \, \, \gamma=0, \, \, \alpha>0.
\end{equation}

Further let \begin{equation}\Lambda_{\nu}(t)=\int_{t}^{1}\frac{\lambda(x)}{x^{1/{\nu}}}dx, \, \, \, \, \, {\nu}>0,\end{equation}
and
\begin{equation}\Pi_{\mu,\nu}(t)=\left\{
                                      \begin{array}{ll}
                                        \int_{t}^{1}\Lambda_{\nu}(x){x}^{1/{\nu}-1-1/{\mu}}dx, & \hbox{$\gamma>0$,} \\
                                        \Lambda_{\alpha}(t), & \hbox{$\gamma=0$, $(\mu=0, \, \nu=\alpha>0)$.}
                                      \end{array}
                                    \right.
\end{equation}
For the function $\Pi_{\mu,\nu}(t)$, we define
\begin{equation}{\mathfrak{M}}_{\Pi_{\mu,\nu}}(h_{\delta})=\left\{
                  \begin{array}{ll}
                    \Re\int_{0}^{1}t^{1/{\mu}-1}{\Pi_{\mu,\nu}}(t)\left[h'_{\delta}(tz)-\frac{(1-\delta)-(1+\delta)t}{(1-\delta)(1+t)^{3}}\right]dt, & \hbox{$\gamma>0$,} \\

                    \Re\int_{0}^{1}t^{1/{\alpha}-1}\Pi_{0,\alpha}(t)\left[{h'_{\delta}(tz)}-\frac{(1-\delta)-(1+\delta)t}{(1-\delta)(1+t)^{3}}\right]dt, & \hbox{$\gamma=0$,}
                  \end{array}
                \right.
\end{equation}
where $h_{\delta}(z)$ is defined as
\begin{equation}h_{\delta}(z)=\frac{z\left(1+\frac{\epsilon+2\delta-1}{2-2\delta}z\right)}{(1-z)^{2}}, \, \, \, \, |\epsilon|=1.\end{equation}

With these notations, we are now in a position to state our first result, which generalizes many earlier results in this direction.\\


\section{Main results}
\setcounter{equation}{0}

\noindent\textbf{Theorem 3.1.} Let $\mu\geq0$, $\nu\geq0$ satisfy (2.1) . Define $\beta<1$ by
\begin{equation}
\frac{\beta-\frac{1}{2}}{(1-\beta)}=-\int_{0}^{1}\lambda(t)q(t)dt,
\end{equation}
where $q(t)$ is the solution of the initial-value problem (2.5). Further for $\Lambda_{\nu}(t)$ and $\Pi_{\mu,\nu}(t)$ defined by (2.8) and (2.9) respectively, assume that $\displaystyle t^{1/{\nu}}\Lambda_{\nu}(t)\rightarrow0$, and $\displaystyle t^{1/{\nu}}\Pi_{\mu,\nu}(t)\rightarrow0$ as $t\rightarrow0^{+}$. Then for $\delta\in[0,\frac{1}{2}]$, $V_{\lambda}({\mathcal{W}}_{\beta}(\alpha, \gamma)) \subset K(\delta)$ if and only if ${\mathfrak{M}}_{\Pi_{\mu,\nu}}(h_{\delta})\geq0$, where ${\mathfrak{M}}_{\Pi_{\mu,\nu}}(h_{\delta})$ and $h_{\delta}$ are defined by equations (2.10) and (2.11) respectively.\\

\noindent\textbf{Proof.} As the case $\gamma=0$ ($\mu=0$, $\nu=\alpha$) corresponds to the Theorem 2.3 in \cite{balasamy}, so we will prove the result only when $\gamma>0$.\\
Let $$H(z)=(1-\alpha+2\gamma)\frac{f(z)}{z}+(\alpha-2\gamma)f'(z)+\gamma zf''(z).$$
Since $\mu+\nu=\alpha-\gamma$ and $\mu\nu=\gamma$, therefore
\begin{eqnarray*}H(z) & = & (1+\gamma-(\alpha-\gamma))\frac{f(z)}{z}+(\alpha-\gamma-\gamma)f'(z)+\gamma zf''(z)\\
& = & (1+\mu\nu-\mu-\nu)\frac{f(z)}{z}+(\mu+\nu-\mu\nu)f'(z)+\mu\nu zf''(z).\end{eqnarray*}
Writing $f(z)=z+\sum_{n=2}^{\infty}a_{n}z^{n}$, we obtain from (2.2)
\begin{equation}H(z)=1+\sum_{n=1}^{\infty}a_{n+1}{(n\nu+1)(n\mu+1)}z^{n}=f'(z)\ast\phi_{\mu,\nu}(z),\end{equation}
and (2.3) gives that
\begin{equation}f'(z)=H(z)\ast\psi_{\mu,\nu}(z).\end{equation}
Now, for $f\in{\mathcal{W}_{\beta}(\alpha, \gamma)}$, we have
$$\Re\left\{e^{i\phi}\frac{H(z)-\beta}{1-\beta}\right\}>0.$$
Thus, in the view of the Theorem 1.1, we may confine oueselves to functions $f\in{\mathcal{W}_{\beta}(\alpha, \gamma)}$ for which
$$H(z)=\beta+(1-\beta)\left(\frac{1+xz}{1+yz}\right), \, \, \, |x|=|y|=1.$$
Thus (3.3) gives
\begin{equation}f'(z)=\left((1-\beta)\frac{1+xz}{1+yz}+\beta\right)\ast\psi_{\mu,\nu}(z),\end{equation}
and therefore
\begin{equation}\frac{f(z)}{z}=\frac{1}{z}\int_{0}^{z}\left((1-\beta)\frac{1+xw}{1+yw}+\beta\right)dw\ast\psi(z).\end{equation}
Here $\psi:=\psi_{\mu,\nu}$.\\

A well-known result from the theory of convolutions [9, Pg 94] (also see \cite{ruscheweyh1}) states that
$$F\in K(\delta) \, \, \, \Leftrightarrow \, \, \,  \frac{1}{z}(zF'\ast h_{\delta})(z)\neq0, \, \, \, z\in E,$$
where
$$h_{\delta}(z)=\frac{z\left(1+\frac{\epsilon+2\delta-1}{2-2\delta}z\right)}{(1-z)^{2}}, \, \, \, \, |\epsilon|=1.$$
Hence $F\in K(\delta)$ if and only if
$$0\neq\frac{1}{z}(V_{\lambda}(f)(z)\ast z{h'}_{\delta}(z))=\frac{1}{z}\left[\int_{0}^{1}\lambda(t)\frac{f(tz)}{t}dt\ast z{h'}_{\delta}(z)\right]=\int_{0}^{1}\frac{\lambda(t)}{1-tz}dt\ast\frac{f(z)}{z}\ast{{h'}_{\delta}(z)} $$
Using (3.5), we have
\begin{eqnarray*}0 & \neq & \int_{0}^{1}\frac{\lambda(t)}{1-tz}dt\ast\left[\frac{1}{z}\int_{0}^{z}\left((1-\beta)\frac{1+xw}{1+yw}+\beta\right)dw\ast\psi(z)\right]\ast{{h'}_{\delta}(z)}\\
& = & \int_{0}^{1}\frac{\lambda(t)}{1-tz}dt\ast{{h'}_{\delta}(z)}\ast\left[\frac{1}{z}\int_{0}^{z}\left((1-\beta)\frac{1+xw}{1+yw}+\beta\right)dw\right]\ast\psi(z)\\
& = & \int_{0}^{1}\lambda(t){{h'}_{\delta}(tz)}dt\ast(1-\beta)\left[\frac{1}{z}\int_{0}^{z}\left(\frac{1+xw}{1+yw}+\frac{\beta}{(1-\beta)}\right)dw\right]\ast\psi(z)\\
& = & (1-\beta)\left[\int_{0}^{1}\lambda(t){{h'}_{\delta}(tz)}dt+\frac{\beta}{(1-\beta)}\right]\ast\frac{1}{z}\int_{0}^{z}\frac{1+xw}{1+yw}dw\ast\psi(z)\\
& = & (1-\beta)\left[\int_{0}^{1}\lambda(t)\left(\frac{1}{z}\int_{0}^{z}{{h'}_{\delta}(tw)}dw\right)dt+\frac{\beta}{(1-\beta)}\right]\ast\frac{1+xz}{1+yz}\ast\psi(z).\end{eqnarray*}
This holds if and only if [11, p. 23]
\begin{eqnarray*} &  & \Re(1-\beta)\left[\int_{0}^{1}\lambda(t)\left(\frac{1}{z}\int_{0}^{z}{{h'}_{\delta}(tw)}dw\right)dt+\frac{\beta}{(1-\beta)}\right]\ast\psi(z)\geq\frac{1}{2},\\ & \Leftrightarrow & \Re(1-\beta)\left[\int_{0}^{1}\lambda(t)\left(\frac{1}{z}\int_{0}^{z}{{h'}_{\delta}(tw)}dw\right)dt+\frac{\beta}{(1-\beta)}-\frac{1}{2(1-\beta)}\right]\ast\psi(z)\geq0,\\
& \Leftrightarrow & \Re\left[\int_{0}^{1}\lambda(t)\left(\frac{1}{z}\int_{0}^{z}{{h'}_{\delta}(tw)}dw\right)dt+\frac{\beta-\frac{1}{2}}{(1-\beta)}\right]\ast\psi(z)\geq0,\\
& \Leftrightarrow &
\Re\left[\int_{0}^{1}\lambda(t)\left(\frac{1}{z}\int_{0}^{z}{{h'}_{\delta}(tw)}dw-q(t)\right)dt\right]\ast\psi(z)\geq0, \, \, \, \, \text{(using (3.1))}
\end{eqnarray*}

\begin{eqnarray*}
& \Leftrightarrow &
\Re\left[\int_{0}^{1}\lambda(t)\left({{h'}_{\delta}(tz)}-q(t)\right)dt\right]\ast\frac{1}{z}\int_{0}^{z}\psi(w)dw\geq0,\\
& \Leftrightarrow &
\Re\left[\int_{0}^{1}\lambda(t)\left({{h'}_{\delta}(tz)}-q(t)\right)dt\right]\ast\sum_{n=0}^{\infty}\frac{z^n}{(n\nu+1)(n\mu+1)}\geq 0, \, \, \, \, \text{(using (2.3))}\\
& \Leftrightarrow &
\Re\int_{0}^{1}\lambda(t)\left(\sum_{n=0}^{\infty}\frac{z^n}{(n\nu+1)(n\mu+1)}\ast{{h'}_{\delta}(tz)}-q(t)\right)dt \geq 0,\\
& \Leftrightarrow &
\Re\int_{0}^{1}\lambda(t)\left(\int_{0}^{1}\int_{0}^{1}\frac{d{\eta}d{\zeta}}{(1-{\eta}^{\nu}{\zeta}^{\mu}z)}\ast{{h'}_{\delta}(tz)}-q(t)\right)dt\geq 0,\\
& \Leftrightarrow &
\Re\int_{0}^{1}\lambda(t)\left(\int_{0}^{1}\int_{0}^{1}{{h'}_{\delta}(tz{\eta}^{\nu}{\zeta}^{\mu})}d{\eta}d{\zeta}-q(t)\right)dt \geq 0,
\end{eqnarray*}
which can also be written as
$$\Re\int_{0}^{1}\lambda(t)\left(\int_{0}^{1}\int_{0}^{1}\frac{1}{\mu\nu}{{h'}_{\delta}(tzuv)}u^{1/{\nu}-1}v^{1/{\mu}-1}dvdu-q(t)\right)dt\geq0.$$
Writing $w=tu$, we get
$$\Re\int_{0}^{1}\frac{\lambda(t)}{t^{1/{\nu}}}\left[\int_{0}^{t}\int_{0}^{1}{{h'}_{\delta}(wzv)}w^{1/{\nu}-1}v^{1/{\mu}-1}dvdw-{\mu\nu}{t^{1/{\nu}}}q(t)\right]dt\geq0.$$
An integration by parts with respect to $t$ and (2.5) gives
$$\Re\int_{0}^{1}\Lambda_{\nu}(t)\left[\int_{0}^{1}{{h'}_{\delta}(tzv)}t^{1/{\nu}-1}v^{1/{\mu}-1}dv-t^{1/{\nu}-1}\int_{0}^{1}\frac{1-\delta-(1+\delta)st}{(1-\delta)(1+st)^3}s^{1/{\mu}-1}ds\right]dt\geq0.$$
Again writing $w=vt$ and $\eta=st$ above inequality reduces to
$$\Re\int_{0}^{1}\Lambda_{\nu}(t)t^{1/{\nu}-1/{\mu}-1}\left[\int_{0}^{t}{{h'}_{\delta}(wz)}w^{1/{\mu}-1}dw-\int_{0}^{t}\frac{1-\delta-(1+\delta){\eta}}{(1-\delta)(1+\eta)^3}{\eta}^{1/{\mu}-1}d\eta\right]dt\geq0,$$
which after integration by parts with respect to $t$ yields
$$\Re\int_{0}^{1}\Pi_{\mu,\nu}(t)t^{1/{\mu}-1}\left[{{h'}_{\delta}(tz)}-\frac{1-\delta-(1+\delta)t}{(1-\delta)(1+t)^3}\right]dt\geq0.$$
Thus $F\in K(\delta)$ if and only if ${\mathfrak{M}}_{\Pi_{\mu,\nu}}(h_{\delta})\geq0$.\\

Finally, to prove the sharpness, let $f\in{\mathcal{W}}_{\beta}(\alpha, \gamma)$ be of the form for which
$$(1-\alpha+2\gamma)\frac{f(z)}{z}+(\alpha-2\gamma)f'(z)+\gamma zf''(z)=\beta+(1-\beta)\frac{1+z}{1-z}.$$
Using a series expansion we obtain that
$$f(z)=z+2(1-\beta)\sum_{n=1}^{\infty}\frac{1}{(n\nu+1)(n\mu+1)}z^{n+1}.$$
Thus
$$F(z)=V_{\lambda}(f)(z)=\int_{0}^{1}\lambda(t)\frac{f(tz)}{t}dt=z+2(1-\beta)\sum_{n=1}^{\infty}\frac{\tau_{n}}{(n\nu+1)(n\mu+1)}z^{n+1},$$
where $\tau_{n}=\int_{0}^{1}\lambda(t)t^{n}dt$. From (2.5), it is a simple exercise to write $q(t)$ in a series expansion as
\begin{equation}
q(t)=1+\frac{1}{1-\delta}\sum_{n=1}^{\infty}\frac{(-1)^{n}(n+1)(n+1-\delta)}{(n\nu+1)(n\mu+1)}t^{n}.
\end{equation}
Now, by (3.1) and (3.6), we have
\begin{eqnarray*}
\frac{\beta-\frac{1}{2}}{1-\beta} & = & -\int_{0}^{1}\lambda(t)q(t)dt\\
 & = & -\int_{0}^{1}\lambda(t)\left[1+\frac{1}{1-\delta}\sum_{n=1}^{\infty}\frac{(-1)^{n}(n+1)(n+1-\delta)}{(n\nu+1)(n\mu+1)}t^{n}\right]dt\\
 & = & -1-\frac{1}{1-\delta}\sum_{n=1}^{\infty}\frac{(-1)^{n}(n+1)(n+1-\delta)}{(n\nu+1)(n\mu+1)}\int_{0}^{1}\lambda(t)t^{n}dt.
\end{eqnarray*}
Therefore
\begin{equation}
\frac{1}{2(1-\beta)}=-\frac{1}{1-\delta}\sum_{n=1}^{\infty}\frac{(-1)^{n}(n+1)(n+1-\delta)\tau_{n}}{(n\nu+1)(n\mu+1)}.
\end{equation}
Finally, we see that
$$F'(z)=1+2(1-\beta)\sum_{n=1}^{\infty}\frac{(n+1)\tau_{n}}{(n\nu+1)(n\mu+1)}z^{n}.$$
Therefore
$$(zF'(z))'=1+2(1-\beta)\sum_{n=1}^{\infty}\frac{(n+1)^{2}\tau_{n}}{(n\nu+1)(n\mu+1)}z^{n}.$$
For $z=-1$, we have
\begin{eqnarray*}
(zF')'(-1) & = & 1+2(1-\beta)\sum_{n=1}^{\infty}\frac{(-1)^{n}(n+1)^2\tau_{n}}{(n\nu+1)(n\mu+1)}\\
 & = & 1+2(1-\beta)\sum_{n=1}^{\infty}\frac{(-1)^{n}(n+1)(n+1-\delta)\tau_{n}}{(n\nu+1)(n\mu+1)}+2(1-\beta)\sum_{n=1}^{\infty}\frac{(-1)^{n}\delta(n+1)\tau_{n}}{(n\nu+1)(n\mu+1)}\\
 & = & 1-(1-\delta)+\delta2(1-\beta)\sum_{n=1}^{\infty}\frac{(-1)^{n}(n+1)\tau_{n}}{(n\nu+1)(n\mu+1)} \, \, \, \, (\text{Using} (3.7))\\
 & = & \delta\left(1+2(1-\beta)\sum_{n=1}^{\infty}\frac{(-1)^{n}(n+1)\tau_{n}}{(n\nu+1)(n\mu+1)}\right)\\
 & = & \delta F'(-1).\end{eqnarray*}
Thus $(zF'(z))'/F'(z)$ at $z=-1$ equals $\delta$. This implies that the result is sharp for the order of convexity.

\section{Consequences of Theorem 3.1 }
\setcounter{equation}{0}

To obtain a sufficient condition for the convexity of order $\delta$ of the integral transform (1.1) by a much easier method, we present the following theorem.\\

\textbf{Theorem 4.1.} Let $\Lambda_{\nu}(t)$, $\Pi_{\mu,\nu}(t)$ be integrable on [0,1] and positive on (0,1). Also, suppose that $\displaystyle t^{1/{\nu}}\Lambda_{\nu}(t)\rightarrow0$, and $\displaystyle t^{1/{\nu}}\Pi_{\mu,\nu}(t)\rightarrow0$ as $t\rightarrow0^{+}$. Assume further that $\mu\geq1$ and
\begin{equation}
\frac{\left(-t{{\Pi}'}_{\mu,\nu}(t)+\left(1-\frac{1}{\mu}\right)\Pi_{\mu,\nu}(t)\right)}{(1+t)(1-t)^{1+2\delta}} \, \, \, \text{is decreasing on (0,1).}
\end{equation}
For $\delta\in[0,1/2]$, if $\beta$ satisfies (3.1), then $V_{\lambda}(f)\in K(\delta)$ for $f\in{\mathcal{W}}_{\beta}(\alpha, \gamma)$.\\

\textbf{Proof.} For $\gamma>0$, integration by parts with respect to $t$ yields
\begin{small}
\begin{eqnarray*}
&   & \int_{0}^{1}t^{\frac{1}{\mu}-1}\Pi_{\mu,\nu}(t)\left(\Re{\left({h'}_{\delta}(tz)\right)}-\frac{1-\delta-(1+\delta)t}{(1-\delta)(1+t)^{3}}\right)dt\\
& = & \int_{0}^{1}t^{\frac{1}{\mu}-1}\Pi_{\mu,\nu}(t)\frac{d}{dt}\left(\Re{\frac{{h}_{\delta}(tz)}{z}}-\frac{t(1-\delta(1+t))}{(1-\delta)(1+t)^{2}}\right)dt\\
& = &
\int_{0}^{1}t^{\frac{1}{\mu}-1}\left(-t{\Pi'}_{\mu,\nu}(t)+\left(1-\frac{1}{\mu}\right)\Pi_{\mu,\nu}(t)\right)\left(\Re{\frac{{h}_{\delta}(tz)}{tz}}-\frac{1-\delta(1+t)}{(1-\delta)(1+t)^{2}}\right)dt.
\end{eqnarray*}
\end{small}
\indent Also for $\mu\geq1$, the function $t^{{1/{\mu}}-1}$ is decreasing on (0,1). Thus, the condition (4.1) along with Theorem 1 from \cite{fournier} yields
$$\int_{0}^{1}t^{\frac{1}{\mu}-1}\Pi_{\mu,\nu}(t)\left(\Re{\left({h'}_{\delta}(tz)\right)}-\frac{1-\delta-(1+\delta)t}{(1-\delta)(1+t)^{3}}\right)dt>0.$$
Thus, an application of Theorem 3.1 evidently leads to the desired result.  \hspace{2.5cm}$\square$\\

Below, we obtain the conditions to ensure convexity of $V_{\lambda}(f)$. As defined in (2.8) and (2.9), for $\gamma>0$,
$$\Pi_{\mu,\nu}(t)=\int_{t}^{1}\Lambda_{\nu}(x){x}^{1/{\nu}-1-1/{\mu}}dx, \, \, \, \text{and} \, \, \, \Lambda_{\nu}(t)=\int_{t}^{1}\frac{\lambda(x)}{x^{1/{\nu}}}dx.$$
In order to apply Theorem 4.1, we have to prove that the function $$k(t)=\frac{\left(t^{\frac{1}{\nu}-\frac{1}{\mu}}{{\Lambda}}_{\nu}(t)+\left(1-\frac{1}{\mu}\right)\Pi_{\mu,\nu}(t)\right)}{(1+t)(1-t)^{1+2\delta}}:=\frac{p(t)}{(1+t)(1-t)^{1+2\delta}}$$ is decreasing in (0,1). Since $k(t)>0$ and
\begin{eqnarray*}\frac{k'(t)}{k(t)} & = & \frac{p'(t)}{p(t)}+\frac{2(t+\delta(1+t))}{1-t^{2}}\\
& = &  \frac{2(t+\delta(1+t))}{(1-t^{2})p(t)}\left[\frac{(1-t^{2})p'(t)}{2(t+\delta(1+t))}+{p(t)}\right]=\frac{2(t+\delta(1+t))}{(1-t^{2})p(t)}\left[q(t)\right] \, \, \, \text{(say)}.\end{eqnarray*}
Thus to prove that $k'(t)\leq0$, it is enough to prove that $q(t)\leq0$. Since $q(1)=0$, so it remains to show that $q(t)$ is increasing over (0,1). Now
$$ q'(t) =\frac{(1+t)}{2(t+\delta(1+t))^2}\left[(1-t)(t+\delta(1+t))p''(t)-{(1-t-\delta(1+t))(1+2\delta)}p'(t)\right].$$
So, $q'(t)\geq0$ for $t\in (0,1)$ is equivalent to the inequality $r(t)\geq0$, where
$$r(t)=(1-t)(t+\delta(1+t))p''(t)-{(1-t-\delta(1+t))(1+2\delta)}p'(t)$$
By using the idea similar to the one used to prove Theorem 3.1 in \cite{bala-dj}, we can write
\begin{small}
\begin{equation}
r(t)=-{\lambda(t)}t^{1-\frac{1}{\mu}}\left[\left(\frac{1}{\nu}-\frac{1}{\mu}-1\right)X(t)+Z(t)+\frac{t{\lambda}'(t)}{\lambda(t)}X(t)\right]+\left[\left(\frac{1}{\nu}-\frac{1}{\mu}-1\right)X(t)+Z(t)\right]\left(\frac{1}{\nu}-1\right)t^{\frac{1}{\nu}-\frac{1}{\mu}-1}{\int_{t}^{1}A(s)ds}
\end{equation}
\end{small}
where,
\begin{eqnarray}
        A(t) &=& \lambda(t)t^{-1/{\nu}}, \nonumber\\
        X(t) &=& (1-t)(t+\delta(1+t)), \nonumber\\
        Z(t) &=& -t(1-t-\delta(1+t))(1+2\delta).
\end{eqnarray}
Clearly, $A(t)>0$ and $X(t)>0$ for all $t\in(0,1).$\\
Thus, $r(t)$ is non-negative if
\begin{small}
\begin{equation}
\left(\frac{1}{\nu}-\frac{1}{\mu}-1\right)X(t)+Z(t)+\frac{t{\lambda}'(t)}{\lambda(t)}X(t)\leq 0 \, \, \, \text{and} \, \, \,  \left[\left(\frac{1}{\nu}-\frac{1}{\mu}-1\right)X(t)+Z(t)\right]\left(\frac{1}{\nu}-1\right)\geq0.
\end{equation}
\end{small}
Since $\nu\geq1$, we can rewrite the condition (4.4) as follows :
\begin{small}
\begin{equation}
\frac{t{\lambda}'(t)}{\lambda(t)}\leq 2+\frac{1}{\mu}-\frac{1}{\nu}-\left(\frac{X(t)+Z(t)}{X(t)}\right) \, \, \, \text{and} \, \, \,  \frac{1}{\nu}-\frac{1}{\mu}-2\leq-\left(\frac{X(t)+Z(t)}{X(t)}\right).
\end{equation}
\end{small}
In view of the fact that $X(t)+Z(t)$ and $X(t)$ are non-negative on (0,1), the above inequality further reduces to
\begin{equation}
\frac{t{\lambda}'(t)}{\lambda(t)}\leq 2+\frac{1}{\mu}-\frac{1}{\nu} \, \, \, \text{and} \, \, \,  \frac{1}{\nu}-\frac{1}{\mu}-2\leq0.
\end{equation}
For $\mu\geq1$, condition (2.1) implies $\nu\geq\mu\geq1$. Thus, condition (4.6) implies that $r(t)$ is non-negative if
\begin{equation}
\frac{t{\lambda}'(t)}{\lambda(t)}\leq 2+\frac{1}{\mu}-\frac{1}{\nu},  \hspace{1cm} \nu\geq\mu\geq1.
\end{equation}
These conditions leads to the following theorem.\\

\textbf{Theorem 4.2.} Assume that both $\Lambda_{\nu}(t)$, $\Pi_{\mu,\nu}(t)$ are integrable on [0,1] and positive on (0,1).  Let $\lambda(t)$ be a non-negative real-valued integrable function on [0,1] and satisfy the condition
\begin{equation}
\frac{t{\lambda}'(t)}{\lambda(t)}\leq 2+\frac{1}{\mu}-\frac{1}{\nu}, \hspace{1cm} \nu\geq\mu\geq1.
\end{equation}
Let $f\in{\mathcal{W}_{\beta}(\alpha, \gamma)}$ and $\beta<1$ with
$$\frac{\beta-\frac{1}{2}}{(1-\beta)}=-\int_{0}^{1}\lambda(t)q(t)dt,$$ where $q(t)$ is defined by (2.6). Then $F(z)=V_{\lambda}(f)(z)\in K(\delta)$ for $\delta\in[0,1/2]$. The conclusion does not hold for smaller values of $\beta$.\\

\noindent On the other hand, when $\gamma=0$  ($\mu=0$, $\nu=\alpha>0$), so we get the following result.\\

\noindent\textbf{Theorem 4.3.} Let $\lambda(t)$ be a non-negative real-valued integrable function on [0,1]. Assume that both $\Lambda_{\alpha}(t)$, $\Pi_{0,\alpha}(t)$ are integrable on [0,1] and positive on (0,1). Let $\lambda(1)=0$ and $\lambda$ satisfies the condition
\begin{equation}
t{\lambda}''(t)-\frac{1}{\alpha}{\lambda}'(t)\geq0, \hspace{3cm} \alpha\geq1.
\end{equation}
Let $f\in{\mathcal{W}_{\beta}(\alpha, \gamma)}$ and $\beta<1$ with
$$\frac{\beta-\frac{1}{2}}{(1-\beta)}=-\int_{0}^{1}\lambda(t)q_{\alpha}(t)dt,$$ where $q_{\alpha}(t)$ is defined by (2.7) with $\delta\in[0,1/2]$. Then $F(z)=V_{\lambda}(f)(z)\in K(\delta)$. The conclusion does not hold for smaller values of $\beta$.\\

\noindent\textbf{Proof.} As in Theorem 3.1, for $\gamma=0$ and $f\in{\mathcal{W}_{\beta}(\alpha, \gamma)}$, we have $V_{\lambda}(f)(z)\in K(\delta)$ if
$$\int_{0}^{1}t^{\frac{1}{\alpha}-1}\Pi_{0,\alpha}(t)\left(\Re{\left({h'}_{\delta}(tz)\right)}-\frac{1-\delta-(1+\delta)t}{(1-\delta)(1+t)^{3}}\right)dt>0,$$
which is equivalent to
$$\int_{0}^{1}t^{\frac{1}{\alpha}-1}\left(t^{1-\frac{1}{\alpha}}{\lambda}(t)+\left(1-\frac{1}{\alpha}\right)\Lambda_{\alpha}(t)
\right)\left(\Re\frac{{h}_{\delta}(tz)}{tz}-\frac{1-\delta(1+t)}{(1-\delta)(1+t)^{2}}\right)dt>0.$$
Since $t^{\frac{1}{\alpha}-1}$ is decreasing on (0,1) for $\alpha\geq1$, thus to apply Theorem 1 in \cite{fournier}, it is enough to show that
$$p(t)=\frac{t^{1-\frac{1}{\alpha}}{\lambda}(t)+\left(1-\frac{1}{\alpha}\right)\Lambda_{\alpha}(t)}{(1+t)(1-t)^{1+2\delta}}:=\frac{k(t)}{(1+t)(1-t)^{1+2\delta}}$$ is decreasing on (0,1). Here, logarithmic differentiation implies that
$$\frac{p'(t)}{p(t)}=\frac{2(t+\delta(1+t))}{(1-t^{2})k(t)}\left[\frac{(1-t^{2})k'(t)}{2(t+\delta(1+t))}+{k(t)}\right].$$
Since $p(t)>0$ for $\alpha\geq1$, thus to prove that $p'(t)\leq0$ on (0,1) it remains to show that
$$r(t)={k(t)}+\frac{(1-t^{2})k'(t)}{2(t+\delta(1+t))}\leq0.$$
Since $r(1)=0$, so $r(t)\leq0$ if $r(t)$ is increasing on (0,1). Thus, $r'(t)$ is non-negative if
$$\frac{t^{\frac{-1}{\alpha}}(1+t)}{2(t+\delta(1+t))}\left\{X(t)t{\lambda}''(t)+\left[\left(1-\frac{1}{\alpha}\right)X(t)+Z(t)\right]{\lambda}'(t)\right\}\geq0,$$
where $X(t)$ and $Z(t)$ are as defined in (4.3). Further simplification yields that
$$t{\lambda}''(t)+\left(\frac{X(t)+Z(t)}{X(t)}-\frac{1}{\alpha}\right){\lambda}'(t)\geq0.$$
Since, ${X(t)}$ and ${X(t)+Z(t)}$ are non-negative in (0,1), thus $r'(t)\geq0$ is equivalent to
$$t{\lambda}''(t)-\frac{1}{\alpha}{\lambda}'(t)\geq0, \hspace{3cm} \alpha\geq1,$$
which completes the proof.\\

\noindent\textbf{Remarks 4.4.} Observe that results in \cite{rmali1} can be obtained from our results by setting $\delta=0$.


\section{Applications}
\setcounter{equation}{0}

In this section, we apply Theorem 4.2 and Theorem 4.3 to obtain certain results regarding convexity of well-known integral operators. The proofs of the following results run on the same lines as given in \cite{rmali1} and hence omitted.\\

Consider $\lambda$ to be defined as
$${\lambda}(t)=(1+c)t^{c}, \hspace{3cm} c>-1.$$
Then the integral transform
\begin{equation}
F_{c}(z)=V_{\lambda}(f)(z)=(1+c)\int_{0}^{1}t^{c-1}f(tz)dt, \hspace{3cm} c>-1,
\end{equation}
is the well-known Bernardi integral operator. The classical Alexander and Libera transforms are special cases of (5.1) with $c=0$ and $c=1$ respectively. For this special case of $\lambda$, the following result holds.

\noindent\textbf{Theorem 5.1.}
Let $c>-1$ and $0<\gamma\leq\alpha\leq1+2\gamma$. Let $\beta<1$ satisfy
$$\frac{\beta-\frac{1}{2}}{1-\beta}=-(1+c)\int_{0}^{1}t^{c}q(t)dt,$$
where $q$ is given by
$$q(t)=\frac{1}{\mu\nu}\int_{0}^{1}\int_{0}^{1}\frac{(1-\delta)-(1+\delta)swt}{(1-\delta)(1+swt)^{3}}s^{1/{\mu}-1}w^{1/{\nu}-1}dsdw.$$
Then for $\delta\in[0,1/2]$, we have $V_{\lambda}({\mathcal{W}}_{\beta}(\alpha, \gamma)) \subset K(\delta)$ provided $c$ satisfies the condition :
\begin{equation}
c\leq 2+\frac{1}{\mu}-\frac{1}{\nu}, \, \, \,  \nu\geq\mu\geq1.
\end{equation}
The value of $\beta$ is sharp.\\

Writing $\alpha=1+2\gamma$, $\gamma>0$ and $\mu=1$ in Theorem 5.1 gives the following criteria of convexity :\\

\noindent\textbf{Corollary 5.2.} Let $-1<c\leq3-1/{\gamma}$ and $\gamma\geq1$. Let $\beta<1$ satisfy
$$\frac{\beta-\frac{1}{2}}{1-\beta}=-(1+c)\int_{0}^{1}t^{c}q(t)dt,$$
where $q$ is given by
$$q(t)=\frac{1}{\mu\nu}\int_{0}^{1}\int_{0}^{1}\frac{(1-\delta)-(1+\delta)swt}{(1-\delta)(1+swt)^{3}}s^{1/{\mu}-1}w^{1/{\nu}-1}dsdw.$$
Then for $\delta\in[0,1/2]$, we have $V_{\lambda}({\mathcal{R}}_{\beta}(\gamma)) \subset K(\delta)$. The value of $\beta$ is sharp.\\

Further, letting $\gamma=1$ and $c=0$ in Corollary 5.2, we have\\

\noindent\textbf{Corollary 5.3.} Let $\beta<1$ satisfy
$$\frac{\beta-\frac{1}{2}}{1-\beta}=\frac{1}{1-\delta}\left(\delta\frac{{\pi}^2}{12}-\log2\right)$$
If $f\in{\mathcal{R}}_{\beta}(1)$, then Alexander transform $\displaystyle F_{0}(z)\equiv A[f](z)=\int_{0}^{1}\frac{f(tz)}{t}dt$ is convex of order $\delta$ where $\delta\in[0,1/2]$. The value of $\beta$ is sharp.\\

\noindent\textbf{Remark 5.4.}
1. For $\delta=0$, $$ \beta_{0}=\frac{1-2\log2}{2-2\log2}=-0.629\ldots.$$
Then, for $f$ satisfying
$$\Re e^{i\phi}\left(f'(z)+zf''(z)-\beta\right)>0, \, z\in E,$$
Alexander transform $A[f]$ is convex. It has been shown in \cite{fournier} that $\beta_{0}$ is the best possible bound here.\\

2. We note that for $\delta=1/2$, $\beta_{1/2}=0.590\ldots$. Then, for $f$ satisfying
$$\Re e^{i\phi}\left(f'(z)+zf''(z)-\beta\right)>0, \, z\in E,$$
Alexander transform $A[f]$ is convex of order $\frac{1}{2}$.\\

While, the case $c=0$ in Theorem 5.1 yields yet another interesting result, which we state as a theorem.\\

\noindent\textbf{Theorem 5.5.}
Let $0<\gamma\leq\alpha\leq1+2\gamma$. If $F\in\mathcal{A}$ satisfies
$$\Re\left(F'(z)+{\alpha}{z}F''(z)+{\gamma}z^{2}F'''(z)\right)>\beta, \hspace{2cm} z\in E,$$
and $\beta<1$ satisfies
$$\frac{\beta-\frac{1}{2}}{1-\beta}=-\int_{0}^{1}q(t)dt,$$
where $q$ is given by
$$q(t)=\frac{1}{\mu\nu}\int_{0}^{1}\int_{0}^{1}\frac{(1-\delta)-(1+\delta)swt}{(1-\delta)(1+swt)^{3}}s^{1/{\mu}-1}w^{1/{\nu}-1}dsdw,$$
then for $\delta\in[0,1/2]$, $F$ belongs to $K(\delta)$. The value of $\beta$ is sharp.\\

To state our next theorem, we define
\begin{equation}{\lambda}(t)=\left\{
    \begin{array}{ll}
      (a+1)(b+1)\frac{t^{a}(1-t^{b-a})}{b-a}, & \hbox{b$\neq$a;} \\
      (a+1)^{2}t^{a}\log(1/t), & \hbox{b=a,}
    \end{array}
  \right.
\end{equation}
where $b>-1$ and $a>-1$.\\
Then,
$$V_{\lambda}(f)(z)=G_{f}(a,b;z)=\left\{
    \begin{array}{ll}
      \frac{(a+1)(b+1)}{b-a}\int_{0}^{1}{t^{a-1}(1-t^{b-a})}f(tz)dt, & \hbox{b$\neq$a;} \\
      {(a+1)^2}\int_{0}^{1}{t^{a-1}}\log(1/t)f(tz)dt, & \hbox{b=a.}
    \end{array}
  \right.
$$

\noindent\textbf{Theorem 5.6.} Let $b>-1$, $a>-1$ and $0<\gamma\leq\alpha\leq1+2\gamma$. Let $\beta<1$ satisfy
$$ \frac{\beta-\frac{1}{2}}{1-\beta}=-\int_{0}^{1}{\lambda}(t)q(t)dt,$$
where $q$ is given by
$$q(t)=\frac{1}{\mu\nu}\int_{0}^{1}\int_{0}^{1}\frac{(1-\delta)-(1+\delta)swt}{(1-\delta)(1+swt)^{3}}s^{1/{\mu}-1}w^{1/{\nu}-1}dsdw.$$
and $\lambda(t)$ is defined by (5.3). If $f\in\mathcal{W}_{\beta}(\alpha, \gamma)$, then the convolution operator $G_{f}(a,b;z)$ belongs to $K(\delta)$ with $\delta\in[0,1/2]$ if
\begin{equation}
a\leq 2+\frac{1}{\mu}-\frac{1}{\nu}, \hspace{3cm} \nu\geq\mu\geq1.
\end{equation}
The value of $\beta$ is sharp.\\

Substituting $\alpha=1+2\gamma$, $\gamma>0$ and $\mu=1$ in Theorem 5.1, gives the following result :\\

\noindent\textbf{Corollary 5.7.} Let $b>-1$, $-1<a\leq3-1/{\gamma}$ and $\gamma\geq1$. Let $\beta<1$ satisfy
$$ \frac{\beta-\frac{1}{2}}{1-\beta}=-\int_{0}^{1}{\lambda}(t)q(t)dt,$$
where $q$ is given by
$$q(t)=\frac{1}{\mu\nu}\int_{0}^{1}\int_{0}^{1}\frac{(1-\delta)-(1+\delta)swt}{(1-\delta)(1+swt)^{3}}s^{1/{\mu}-1}w^{1/{\nu}-1}dsdw.$$
and $\lambda(t)$ is defined by (5.3). If $f\in\mathcal{R}_{\beta}(\gamma)$, then the convolution operator
$G_{f}(a,b;z)$ belongs to $K(\delta)$ with $\delta\in[0,1/2]$. The value of $\beta$ is sharp.\\

While for $\gamma=0$, with an application of Theorem 4.3, we get the following result :\\

\noindent\textbf{Theorem 5.8.} Let $b>-1$, $a>-1$ and $\alpha\geq1$. Let $\beta<1$ satisfy
$$ \frac{\beta-\frac{1}{2}}{1-\beta}=-\int_{0}^{1}{\lambda}(t)q_{\alpha}(t)dt,$$
where $q_{\alpha}$ is given by
$$q_{\alpha}(t)=\frac{1}{\alpha}\int_{0}^{1}\frac{(1-\delta)-(1+\delta)st}{(1-\delta)(1+st)^{3}}s^{1/{\alpha}-1}ds$$
and $\lambda(t)$ is defined by (5.3). If $f\in\mathcal{P}_{\beta}(\alpha)$, then the convolution operator
$G_{f}(a,b;z)$ belongs to $K(\delta)$ with $\delta\in[0,1/2]$ if one of the following conditions holds :\\
(i) $-1<a\leq0$ and $a=b$, or\\
(ii) $-1<a\leq0$ and $-1<a<b\leq1+1/{\alpha}$.\\
The value of $\beta$ is sharp.\\

Now, we define
$${\lambda}(t)=\frac{(1+a)^{p}}{\Gamma(p)}t^{a}\left(\log(1/t)\right)^{p-1}, \, \, \, a>-1, \, \, p\geq0.$$
In this case, $V_{\lambda}$ reduces to the Komatu operator \cite{komatu}
$$V_{\lambda}(f)(z)=\frac{(1+a)^{p}}{\Gamma(p)}\int_{0}^{1}\left(\log\left(\frac{1}{t}\right)\right)^{p-1}t^{a-1}f(tz)dt, \, \, \, a>-1, \, \, p\geq0.$$
For $p = 1$ Komatu operator gives the Bernardi integral operator. For this $\lambda$, the following result holds.

\noindent\textbf{Theorem 5.9.}
Let $a>p-2>-1$ and $0<\gamma\leq\alpha\leq1+2\gamma$. Let $\beta<1$ satisfy
$$ \frac{\beta-\frac{1}{2}}{1-\beta}=-\int_{0}^{1}{\lambda}(t)q(t)dt,$$
where $q$ is given by
$$q(t)=\frac{1}{\mu\nu}\int_{0}^{1}\int_{0}^{1}\frac{(1-\delta)-(1+\delta)swt}{(1-\delta)(1+swt)^{3}}s^{1/{\mu}-1}w^{1/{\nu}-1}dsdw.$$
If $f\in\mathcal{W}_{\beta}(\alpha, \gamma)$, then the function
$${\Phi}_{p}(a;z)\ast f(z)=\frac{(1+a)^{p}}{\Gamma(p)}\int_{0}^{1}\left(\log\left(\frac{1}{t}\right)\right)^{p-1}t^{a-1}f(tz)dt$$
belongs to $K(\delta)$ with $\delta\in[0,1/2]$ if
\begin{equation}
a\leq2+\frac{1}{\mu}-\frac{1}{\nu}, \hspace{3cm} \nu\geq\mu\geq1.
\end{equation}
The value of $\beta$ is sharp.\\


\end{document}